\documentclass[11pt]{article}   % For Latex2e
\usepackage{amsthm, amsmath,amsfonts, color}   % For Latex2e
\theoremstyle{plain}
\newtheorem{theorem}{Theorem}[section]
\newtheorem{lemma}[theorem]{Lemma}
\newtheorem{proposition}[theorem]{Proposition}

\theoremstyle{definition}
\newtheorem{example}[theorem]{Example}

\theoremstyle{remark}
\newtheorem{remark}[theorem]{Remark}

\newenvironment{ack}{\noindent{\bf Acknowledgements}.}{}

\def\KK{{\mathbb K}}

\def\AA{{\mathcal A}}

\def\BB{{\mathcal B}}

\def\t{{\mathbf t}}
\def\e{{\mathbf e_{\AA_0,\ldots,\AA_k}}}
\def\eB{{\mathbf e_{\BB_0,\ldots,\BB_s}}}
\def\d{{\mathbf d_{\AA_0,\ldots,\AA_k}}}

\def\X{{\mathcal X}}

\def\CC{{\mathbb C}}

\def\QQ{{\mathbb Q}}
\def\RR{{\mathbb R}}
\def\ZZ{{\mathbb Z}}

\def\Ch{{\mathcal C}h}

\title{Rational Formulas for Traces in zero-dimensional Algebras}
\author{
Carlos D'Andrea\thanks{Supported by the Programa Ram{\'o}n y Cajal,
Ministerio de Educaci{\'o}n y Ciencia, Spain and also by the
Research Project MTM2007-67493 of the Ministerio de Educaci\'on y
Ciencia, Spain. }  \and Gabriela Jeronimo\thanks{Partially supported
by the Argentinian research grants UBACyT X847 (2006-2009), CONICET
PIP 5852/05 and ANPCYT PICT 2005 17-33018.} \footnote{Corresponding
author.}}
\date{}

\begin{document}
\maketitle

\begin{abstract}
We present a rational expression for the trace of the multiplication
map $\mbox{Times}_r:A\to A$ in a finite-dimensional algebra $A:=
\KK[x_1,\ldots,x_n] /\mathcal{I}$ in terms of the generalized Chow
form of $\mathcal{I}$. Here, $\mathcal{I}\subset
\KK[x_1,\ldots,x_n]$ is a zero-dimensional ideal, $\KK$ is a field
of characteristic zero, and $r(x_1,\dots, x_n)$ a rational function
whose denominator is not a zero divisor in $A$.  If $\mathcal{I}$ is
a complete intersection in the torus, we get numerator and
denominator formulas for traces in terms of sparse resultants.
\end{abstract}

\section{Introduction}
Traces in finite dimensional algebras play a fundamental role in
Commutative Algebra and Algebraic Geometry. Recent applications of
traces include the evaluation of symmetric functions, the effective
Nullstellensatz, the computation of radicals of ideals and
algorithms for solving polynomial systems \cite{ABRW,AS,
BW,DG,FGS,JRS,KP,SS,Kun}. For developments in the analytical
counterpart of this algebraic tool (residues) see
\cite{AGV,AY,BGVY,BY1,BY2,BY3,CDS1,CDS,CM,elk,Tsi}.

\par The main results of this paper are explicit rational formulas
for the computation of traces of the multiplication map
$\mbox{Times}_r:A\to A$ in terms of the generalized Chow form of
$\mathcal{I}$ in the case $\mathcal{I}\subset \KK[x_1,\ldots,x_n]$ is a zero-dimensional ideal and $r(x_1,\dots, x_n)$
a rational function whose denominator is not a
zero divisor in $A:= \KK[x_1,\ldots,x_n] /\mathcal{I}$.  If $\mathcal{I}$ is a complete intersection
in the torus, we get numerator and denominator formulas for traces
in terms of \textit{sparse resultants}.
\par The importance of having such formulas is that they allow to develop computational techniques for
solving polynomial equations  as in \cite{JRS}, to bound arithmetic
aspects of the membership problem as in \cite{elk} or to compute
invariants in ``general'' (as opposed to ``generic'') situations
like the computations given in \cite{ped,EY} to mention some
examples. In all these cases, from generic formulae like the results
presented here, one can perform a suitable specialization and get
results in a general case. 
\par\smallskip To state properly our results, let $\KK$ be a field of characteristic zero,  $S:=\KK[x_1,\ldots,x_n],$
$\mathcal{I}\subset S$ a zero-dimensional ideal, and $A:=
S/\mathcal{I}.$ Let $p,q\in S$ such that $q$ is not a zero
divisor in $A.$ Then, $r:=\frac{p}{q}$ induces a
$\KK$-linear map
\begin{equation}\label{uno}
\begin{array}{cccc}
\mbox{Times}_r:&A&\to&A \\
& a&\mapsto&r\cdot a.
\end{array}
\end{equation}
We are interested in the computation of the trace of this map.
First, we relate ${\rm Trace}(\mbox{Times}_r)$ with an algebraic
object depending on $\mathcal{I}$ and $r$: the \textit{generalized
Chow form} of the ideal $\mathcal{I}$ (see \cite{Phi} for the
definition of this eliminating polynomial).
\par Then, we focus on the case where $\mathcal{I}$ is given
by a complete intersection in the torus. We show that, in this
case, the generalized Chow form may be replaced with a sparse
resultant (in the sense of \cite{CLO,GKZ}), and we exhibit numerator and
denominator formulas for traces in terms of sparse resultants,
similar to the denominator formulas obtained for residues in
\cite{CDS}. In particular, our trace formulas can be used for
computing global residues.
\par Finally, we compare our formulas for the denominator of the trace with those
obtained from the formulas proposed in \cite{CDS} for denominators
of residues in the torus.
\par Effective procedures for the computation of traces can be derived from our formulas.
In order to do this, one has to deal with Chow forms and
resultants. In the case $\mathcal{I}$ is a radical ideal, the
generalized Chow form of $\mathcal{I}$ coincides with the
generalized Chow form of the variety $V(\mathcal{I})$ and there
are effective algorithms for its computation (see, for instance,
\cite{JKSS}). Algorithms for the computation of sparse resultants
can be found in \cite{CE,CLO,EP,JS} and the references given therein.

\par On the other hand, our factorization formulas for denominators of
traces and residues could lead to a better understanding of
non-generic situations, namely, those coefficient vectors for which
the denominators vanish.

\par The paper is organized as follows: In Section \ref{2} we show the general formula
for the computation of the trace based on the generalized Chow
form, and then we focus on the case where we have a
generic complete intersection in the torus. We give rational
expressions for both the numerator and the denominator of the
trace. We compare our results with those obtained in \cite{CDS}
for the computation of residues in Section \ref{4}. Section
\ref{5} is concerned with multidimensional residues on affine
space. By using results of Jouanolou we recover known denominator
formulas for residues given in \cite{Tsi,elk,CDS1,CDS} and also give
an algebraic proof of the Euler-Jacobi formula.

\bigskip

\begin{ack}
We are grateful to Jean Pierre Jouanolou for helpful conversations
about his work on discriminants and also to Abdallah Al-Amrani for
sharing with us the notes of the course \cite{jou}. All our examples were computed with the
aid of the software {\tt Maple}.
\end{ack}

\section{Rational expressions for traces}\label{2}
\subsection{A Chow form--based formula}

As in the introduction, let $\KK$ be a field of characteristic
zero, and $S:=\KK[x_1,\ldots,x_n].$ Let $\mathcal{I}$ be a
zero-dimensional ideal of $S,$ and set $A:= S/\mathcal{I}.$
\par Let $p,q\in S$ such that $q$ is not a zero divisor in $A.$
Set $d:=\max\{\deg(p),\deg(q)\},$ where $\deg$ denotes total
degree in the variables $x_1,\ldots,x_n.$ We introduce new
variables $U_\alpha$ for $\alpha\in(\ZZ_{\ge 0})^n,$ $|\alpha|\le
d.$ Let $U:=(U_\alpha,|\alpha|\leq d)$ and
$\Ch_{d,\mathcal{I}}(U)$ be the generalized Chow form of
$\mathcal{I}$ (see \cite{Phi}): if ${\mathcal
U}(x):=\sum_{|\alpha|\leq d}U_\alpha x^\alpha,$ then $
\Ch_{d,\mathcal{I}}(U)=\prod_{\xi\in V(\mathcal{I})}{\mathcal
U}(\xi)^{m(\xi)}, $ where $m(\xi)$ is the multiplicity of $\xi$
with respect to $\mathcal{I},$ i.e. the dimension of the local
ring $A_\xi:=S_\xi / \mathcal{I}S_\xi.$

\begin{theorem}\label{mt}
Let ${\rm Times}_r$ be the map defined in (\ref{uno}). Then,
\begin{equation}\label{doss}
{\rm Trace}({\rm Times}_r)=\frac{\sum_{|\alpha|\leq d}{\,
p_\alpha\,\frac{\partial \Ch_{d,\mathcal{I}}}{\partial
U_\alpha}(q)}}{\Ch_{d,\mathcal{I}}(q)}
\end{equation}
\end{theorem}
Observe that as $q$ is not a zero divisor in $A,$ then the
denominator of (\ref{doss}) is not zero.

%\smallskip
\begin{proof}
Let $T$ be a new variable and set
$$\X(T):=\Ch_{d,\mathcal{I}}(q+T\,p)=C_0+C_1T+\mbox{higher order terms in}\ T.$$
Using the identity $ \Ch_{d,\mathcal{I}}(q+T\,p)=\prod_{\xi\in
V(\mathcal{I})}(q(\xi)+T\,p(\xi))^{m(\xi)},$ we get
\begin{eqnarray*}
\X(T) & =&  \prod_{\xi\in V(\mathcal{I})} q(\xi)^{m(\xi)} +
\Big(\sum_{\xi \in V(\mathcal{I})} m(\xi) p(\xi)
q(\xi)^{m(\xi)-1}\prod_{\xi'\in V(\mathcal{I}), \, \xi'\ne \xi}
q(\xi')^{m(\xi')}\Big) \,  T + {} \\
& & {}+\mbox{higher order terms in}\ T.
\end{eqnarray*}
{}From these identities, recalling that the trace of
$\mbox{Times}_r$ is equal to $\sum_{\xi\in V(I)} m(\xi) r(\xi)$, it
is easy to see that $${\rm
Trace}(\mbox{Times}_r)=\frac{C_1}{C_0}=\frac{\X'(0)}{\X(0)}.$$ By
applying the chain rule to compute the numerator, we get Identity
(\ref{doss}).
\end{proof}

\subsection{Traces in the torus}\label{3} Now we turn our attention
to generic complete intersections in the torus. We will give
rational formulas for traces of linear maps in terms of sparse
resultants.

%\subsection{Rational formulas for traces in the torus}

\medskip

Consider a system of $k$ generic Laurent polynomials $
f_i:=\sum_{a\in \AA_i} c_{ia} \, \mathbf{t}^a$, where
$i=1,\ldots,k,$ $\AA_i$ is a finite subset of $\ZZ^k$, the
coefficients $c_{ia}$ are indeterminates over $\QQ,$
$a=(a_{1},\ldots,a_{k})\in\ZZ^k$ and
$\t^{a}={t_1}^{a_{1}}{t_2}^{a_{2}}\ldots {t_k}^{a_{k}}.$ For
$i=1,\dots, k$, the support of $f_i$ is the set of exponent vectors
$\AA_i\subset\ZZ^k,$ and its Newton polytope is the convex hull
$P_i=conv(\AA_i)\subset \RR^k.$ We assume that the lattice affinely
generated by $\AA_1,\ldots,\AA_k$ (that is, the
$\mathbb{Z}$-submodule of $\mathbb{Z}^k$ generated by the vectors
that are differences of two points in a set $\AA_i$) is a
$k$-dimensional affine sublattice of $\ZZ^k.$
\par Let $K$ be an algebraic
closure of $\QQ(c_{ia})_{1\le i \le k,\, a\in \AA_i}$ and
$S:=K[t_1,{t_1}^{-1},\ldots,t_k,{t_k}^{-1}].$ Then $A:=S/\langle
f_1,\ldots,f_k\rangle$ is a finite-dimensional $K$-vector space
(see \cite{PS}). Let $\AA,\AA'$ be finite subsets of $\ZZ^k$ and
consider a generic rational function of the form
$r(\t):=\frac{p(\t)}{q(\t)}$ where
$p(\t)=\sum_{a\in\AA}p_a\t^a,\,q(\t)=\sum_{a\in\AA'}q_a\t^a,$ and
$p_a,q_a$ are new indeterminates. As $q(\t)$ is a generic
denominator, it is invertible in $A.$ We can then consider the
$K$-linear map $\mbox{Times}_r:A\to A$ as defined in (\ref{uno}).
\par
For any family of finite subsets $\mathcal{B}_1,\dots,
\mathcal{B}_s\subset \ZZ^k$ we denote $L(\mathcal{B}_1,\dots,
\mathcal{B}_s)$ the affine lattice generated by
$\mathcal{B}_1,\dots, \mathcal{B}_s$:
$$L(\mathcal{B}_1,\dots, \mathcal{B}_s)=
\Big\{ \sum_{1\le i\le s} \lambda_i b^{(i)} \mid b^{(i)}\in
\mathcal{B}_i,\, \lambda_i\in \ZZ \hbox{ for all } 1\le i\le s
\hbox{ and } \sum_{1\le i\le s} \lambda_i = 1\Big\}.$$

Let $\AA_0:=\AA\cup\AA'$ and $f_0=\sum_{a\in\AA_0}c_{0a}\t^a,$ where
$c_{0a}$ are new indeterminates.
\par
As in \cite[Section 1]{stu}, for $I\subset\{0,1,\dots,k\}$, the
collection of supports $\{\AA_i\}_{i\in I}$ is said to be
\textit{essential} if $\mbox{rank}(L(\AA_i, i\in I))=\#I-1$ and
$\mbox{rank}(L(\AA_j, j\in J))\geq\#J$ for each proper subset $J$ of
$I$.
\par
Let $M$ be the \textit{mixed volume} of the sequence of polytopes
$P_1,\ldots,P_k$ (see \cite{CLO} for a definition). {}From now on,
we will assume that $M>0$; otherwise, our problem has no interest.
Then, $\AA_0,\AA_1,\ldots,\AA_k$ has a unique essential subset
containing $\AA_0.$ Moreover, the sparse resultant operator ${\rm
Res}_{\AA_0,\ldots,\AA_k}$ as defined in \cite{stu} is not
constantly one and, if $\{\AA_i \}_{i\in I}$ is the essential
subset, it coincides with the resultant ${\rm Res}_{\AA_i, \,i \in
I}$ considered with respect to the lattice $L(\AA_i, i\in I)$ (see
\cite[Corollary 1.1]{stu}).
\par
\medskip
Throughout this section, we are going to use some results and
notation from \cite{min} that we now recall.
\par
For a family of finite sets $\BB_0,\ldots,\BB_s\subset\ZZ^k,$ if
$\{\BB_i\}_{i\in I}$ is the unique essential subfamily, consider the
orthogonal decomposition $L(\BB_0,\ldots,\BB_s)=sat\left(L(\BB_i,\,
i\in I)\right)\oplus L^\perp$, where $sat()$ denotes saturation with
respect to the ambient lattice and $L^\perp$ is the orthogonal
complement (recall that, if $L$ is a sublattice of a lattice
$\mathcal{L}$, the \emph{saturation} of $L$ with respect to
$\mathcal{L}$ is defined as $sat(L) = (\QQ \otimes_\ZZ L) \cap
\mathcal{L}$). If $L^\perp=0,$ we define $\eB:=1.$ Otherwise, denote
with $\pi$ the projection onto the second factor and define $\eB$ as
the normalized mixed volume of the family
$\{\pi(conv(\mathcal{B}_i))\}_{i\notin I}$ in $L^\perp$
(\cite{min}).
%\par
%Let $\d:=[\ZZ^k:L(\AA_0,\ldots,\AA_k)]\,\e.$

\begin{remark}
The reader should be cautious when comparing our statements with
Minimair's results in \cite{min}. Indeed, in Remark $3$ of
\cite{min} the resultant is defined as a power of an irreducible
polynomial. In this paper, a resultant is always an irreducible
polynomial.
\end{remark}

Assuming that $\AA_0,\AA_1,\dots, \AA_k$ has a unique essential
subset containing $\AA_0$, the following Poisson-type product
formula holds (\cite[Lemma $13$]{min}):
\begin{equation}\label{poissontype}
{{\rm Res}_{\AA_0,\ldots,\AA_k}(f_0,f_1,\ldots,f_k)}^\e=\ {\rm
C}\prod_{\beta\in V^*(f_1, \ldots,f_k)}f_0(\beta),
\end{equation}
where ${\rm C}\in K^*,$ and $V^*$ is the set of zeros over $K^*$
with respect to the lattice $L(\AA_0,\ldots,\AA_k).$

For every $\omega \in \ZZ^k$, every set $\AA\subset \ZZ^k$ and any
polynomial $f$ in $k$ variables with support $\AA$, we denote
$a_\AA(\omega) = - \min \{ \langle \omega, v\rangle : v \in \AA\}$,
$\AA^\omega$ the set of points in $\AA$ that lie in the face with
inward normal vector $\omega$, $f^\omega$ the polynomial formed by
the monomials of $f$ lying in $\AA^\omega$, and $H^\omega$ the
lattice of integer points contained in the hyperplane orthogonal to
$\omega$ in $\ZZ^k$.

Under the previous assumptions and notation, if $\AA'_0\subset
\AA_0$ and $f'_0$ is a generic polynomial with support $\AA'_0$, the
main result in \cite{min} establishes a relation between the
specialized resultant ${\rm Res}_{\AA_0, \AA_1,\dots,
\AA_k}(f'_0,f_1,\dots, f_k)$ and ${\rm Res}_{\AA'_0, \AA_1,\dots,
\AA_k}(f'_0,f_1,\dots, f_k)$ (\cite[Theorem 1]{min}):
\begin{equation}\label{mainmin}
\begin{array}{c}
{\rm Res}_{\AA_0, \AA_1,\dots, \AA_k}(f'_0,f_1,\dots, f_k)^{
\mathbf{e}_{\AA_0, \AA_1,\dots, \AA_k}} = {}\hspace{6cm} \\[2mm] {} = {\rm
Res}_{\AA'_0, \AA_1,\dots, \AA_k}(f'_0,f_1,\dots,
f_k)^{\mathbf{e}_{\AA'_0, \AA_1,\dots, \AA_k} [L(\AA_0, \AA_1,\dots,
\AA_k): L(\AA'_0,
\AA_1,\dots, \AA_k)]}\cdot {}\\
 {} \cdot \prod_\omega {\rm Res}_{\AA_1^\omega,\dots,
\AA^\omega_k} (f_1^\omega,\dots,
f_k^\omega)^{\mathbf{e}_{\AA_1^\omega,\dots,
\AA_k^\omega}(a_{\AA_0}(\omega) - a_{\AA'_0}(\omega))
\frac{[H^\omega : L(\AA_1^\omega, \dots, \AA_k^\omega)]}{[\ZZ^k :
L(\AA_0, \AA_1,\dots, \AA_k)]} }, \end{array}
\end{equation}
where the product ranges over all the primitive inward normal
vectors to the facets of the convex hull of $\AA_1 + \cdots +\AA_k$.

\bigskip
Identity (\ref{poissontype}) enables us to use ${\rm
Res}_{\AA_0,\ldots,\AA_k}$ instead of the generalized Chow form used
in the previous section for trace computations. The analogue of
Theorem \ref{mt} is the following:

\begin{theorem}

Under the previous assumptions and notations,
\begin{equation}\label{mtsp}
{\rm Trace}({\rm
Times}_r)=\d\frac{\sum_{a\in\AA}{p_a\,\frac{\partial {\rm
Res}_{\AA_0,\ldots,\AA_k}} {\partial
c_{0a}}(q,f_1,\ldots,f_k)}}{{\rm
Res}_{\AA_0,\ldots,\AA_k}(q,f_1,\ldots,f_k)},
\end{equation}
where $\d:=[\ZZ^k:L(\AA_0,\ldots,\AA_k)]\,\e.$
\end{theorem}

\begin{proof}
Suppose first that $[\ZZ^k:L(\AA_0,\ldots,\AA_k)]=1.$ In this case,
Identity (\ref{poissontype}) reads as follows:
$$
{{\rm Res}_{\AA_0,\ldots,\AA_k}(f_0,f_1,\ldots,f_k)}^\e = \ {\rm
C}\prod_{\gamma\in V(f_1, \ldots,f_k)}f_0(\gamma),
$$
with ${\rm C}\in K^*.$ Now, as in the proof of Theorem \ref{mt}, we
make the substitution $f_0\mapsto q+Tp,$ and get
$$\begin{array}{l}
{{\rm Res}_{\AA_0,\ldots,\AA_k}(q+Tp,f_1,\ldots,f_k)}^\e =
{\rm Res}_{\AA_0,\ldots,\AA_k}(q,f_1,\ldots,f_k)^\e +{}\\[2mm]
\quad {} +\e \,{\rm
Res}_{\AA_0,\ldots,\AA_k}(q,f_1,\ldots,f_k)^{\e-1}\,
\displaystyle\frac{\partial {\rm
Res}_{\AA_0,\ldots,\AA_k}(q+Tp,f_1,\ldots,f_k)}{\partial T}|_{T=0}\,
T + {}
\\[3mm] \quad {} +\mbox{higher order terms in } T.
\end{array}$$
By using the chain rule, we have that $$\frac{\partial {\rm
Res}_{\AA_0,\ldots,\AA_k}(q+Tp,f_1,\ldots,f_k)}{\partial T} =
\sum_{a\in \AA} \frac{\partial {\rm
Res}_{\AA_0,\ldots,\AA_k}}{\partial c_{0a}}(q+Tp,f_1,\ldots,f_k)\,
p_a.$$
\par
The claim follows by substituting this identity in the previous one
and noticing that ${\rm Trace}({\rm Times}_r)$ is the quotient of
the coefficient of the degree $1$ term in $T$ by the coefficient of
degree $0$ in $T$ in the polynomial obtained.
\par
In the general case, we also use Identity (\ref{poissontype}). By
raising both sides of this equality to the power $[\ZZ^k : L(\AA_0,
\dots, \AA_k)]$, the claim is proved as in the previous case if we
show that
\begin{equation}\label{mm}
{\prod_{\beta\in V^*(f_1,
\ldots,f_k)}f_0(\beta)}^{[\ZZ^k:L(\AA_0,\ldots,\AA_k)]}=\prod_{\gamma\in V(f_1,\ldots,f_k)}f_0(\gamma).
\end{equation}
By using normal Smith form reduction (this change preserves
resultants, see \cite{min}), we can suppose that
$L(\AA_0,\ldots,\AA_k)=d_1\ZZ\oplus\ldots\oplus d_k\ZZ,$ with
$[\ZZ^k:L(\AA_0,\ldots,\AA_k)]=d_1\,d_2\ldots d_k.$ Then,
(\ref{mm}) follows straightforwardly from \cite[Corollary
$5$]{min}.
\end{proof}

\begin{example}
{\small
Consider the following trivariate system:
$$\left\{\begin{array}{lll}
f_1&=&c_{11}+c_{12}t_1^2+c_{13}t_2^2\\
f_2&=&c_{21}+c_{22}t_2^2+c_{23}t_3^2\\
f_3&=&c_{31}+c_{32}t_1^2+c_{33}t_2^2+c_{34}t_3^2,
\end{array}
\right.$$ and let $\AA=\{(0,2,0)\},\,\AA':=\{(0,0,0)\}.$ In this
case, the family $\{\AA_0,\AA_1,\AA_2,\AA_3\}$ is essential, so
${\mathbf e}_{ \AA_0,\ldots,\AA_3}=1$, but $[\ZZ^3: L(
{\AA_0,\ldots,\AA_3})]=8.$ A straightforward computation shows
that
$${\rm Trace}(\mbox{Times}_{t_2^2})=8\frac{c_{11}c_{32}c_{23}+c_{21}c_{12}c_{34}-c_{31}c_{12}c_{23}}
{-c_{12}c_{22}c_{34}+c_{12}c_{23}c_{33}-c_{32}c_{13}c_{23}}.
$$
This is a generalization of the example which appears at the end
of \cite{ped}, where the intersection of two perpendicular
cylinders with a sphere is considered:
$$\left\{\begin{array}{lll}
f_1&=&-1+t_1^2+t_2^2\\
f_2&=&-1+t_2^2+t_3^2\\
f_3&=&-1+t_1^2+t_2^2+t_3^2,
\end{array}
\right.$$
and the trace of the multiplication by $t_2^2$ is shown to be $8.$
}
\end{example}

%\subsection{The denominator of the trace in the torus}

In the remaining part of this section we will give an explicit
factorization of the denominator of the trace. The following result
is straightforward due to the irreducibility of the resultant of a
system of generic polynomials (\cite{GKZ}), and the fact that the
degree of the numerator of the right hand side of (\ref{mtsp}) with
respect to the variables $c_{ia},q_a$ is strictly less than the
degree of the denominator.

\begin{lemma}
Suppose that $\AA_0=\AA'$ (i.e. $\AA\subset\AA'$). Then either
the trace is identically zero, or the right hand side of
(\ref{mtsp}) is the irreducible representation of the trace as a
rational function in $\QQ(c_{ia}, p_a,q_a).$
\end{lemma}

Suppose now that $\AA$ is not contained in $\AA'.$ As the trace
is linear in the monomial expansion of $p,$ it is enough to
consider the following situation: $p(\t)=\t^a,$ with
$a\notin\AA'.$
\par
Without loss of generality, let $\{\AA_0,\ldots,\AA_j\}$ $(j\le k)$
be the unique essential subfamily of $\AA_0,\ldots,\AA_k$. Let
$\delta_{\AA'}:={\mathbf
e}_{\AA',\AA_1,\ldots,\AA_j}\mathbf{e}_{\AA_0,\AA_1,\dots,
\AA_j}^{-1}[L(\AA_0,\AA_1,\ldots,\AA_j):L(\AA',\AA_1,\ldots,\AA_j)].$

For each facet of the Minkowski sum $P_1+\ldots+P_j,$  we consider
its primitive inward normal vector $\omega$ and define $\mu_\omega:=
\min\{\langle b,\omega\rangle,\,b\in\AA'\} - \min\{\langle
b,\omega\rangle,\,b\in\AA_0\}. $ Observe that $\mu_\omega\geq0$ and
equality may hold. Set $\delta_{\omega}:=\mu_\omega\, {\mathbf
e}_{{\AA^\omega_1},\ldots,{\AA^\omega_j}}
\mathbf{e}_{\AA_0,\AA_1,\dots, \AA_j}^{-1}
[\omega^{\perp}:L({\AA^\omega_1},\ldots,{\AA^\omega_j})],$ where
 $\omega^\perp$ is the
lattice of integer points contained in the hyperplane orthogonal to
$\omega$ in $L(\AA_0,\AA_1,\ldots,\AA_j).$ With this notation,
Identity (\ref{mainmin}) gives us the following:
\begin{proposition}\label{prp}
In the situation described above, we have that the denominator of (\ref{mtsp}) has the following irreducible factorization:
$${\rm Res}_{\AA',\AA_1,\ldots,\AA_k}(q,f_1,\ldots,f_k)^{\delta_{\AA'}}
\prod_{\omega}{{\rm
Res}_{{\AA^\omega_1},\ldots,{\AA^\omega_k}}({f^\omega_1},\ldots,{f^\omega_k})}^{\delta_\omega}$$
where $\omega$ ranges over the primitive inward normal vectors of
the facets of $P_1+\ldots+P_j$.
\end{proposition}

\begin{example}
{\small
Consider the following system
\begin{equation}\label{cds}
\left\{\begin{array}{lll}
f_1&=&c_{11}t_1+c_{12}t_1t_2+c_{13}t_2^2 \\
f_2&=&c_{21}t_2+c_{22}t_1t_2+c_{23}t_1^2.
\end{array}\right.
\end{equation}
%We set $q:=q_{(0,0)}+q_{(1,0)}t_1+q_{(0,1)}t_2,\,a:=(2,0).$
We set $q:=q_1+q_2 t_1+q_3 t_2,\,a:=(2,0).$ The Newton polygon
$P_1+P_2$ is a pentagon whose vertices are
$(0,3),(1,3),(3,1),(3,0),(1,1).$ The inward normal vectors $\omega$
of this polygon satisfying $\mu_\omega>0$ are $(-1,-1)$ and
$(-1,0).$ The facet resultants associated with these edges are
$c_{22}c_{12}-c_{13}c_{23}$ and $c_{23}$ respectively and
\begin{equation}\label{cc}
{\rm Res}_{\AA_0,\AA_1,\AA_2}(q,f_1,f_2)={\rm
Res}_{\AA',\AA_1,\AA_2}(q,f_1,f_2)\,
c_{23}\,(c_{22}c_{12}-c_{13}c_{23}).
\end{equation}
This is the irreducible decomposition of the denominator of ${\rm
Trace}(\mbox{Times}_{{\t^a}/{q}}).$ Indeed, computing explicitly, we
get that its numerator is {\small
$$\begin{array}{l}
2c_{13}c_{11}^2c_{22}^3q_2c_{21}q_1-3c_{13}^2c_{11}c_{21}^2q_2c_{23}^2q_1
-2c_{13}c_{11}c_{21}c_{22}^2q_1^2c_{23}c_{12}+2c_{11}^2c_{21}^2q_3^2c_{23}c_{22}c_{12}
-c_{13}c_{11}^2c_{22}^4q_1^2\\
-c_{13}c_{11}^2c_{22}^2q_2^2c_{21}^2-c_{13}c_{11}^2q_3c_{23}q_2c_{22}c_{21}^2
-c_{13}c_{11}^2q_3c_{23}c_{22}^2q_1c_{21}
-c_{11}^2q_3c_{22}^3c_{12}q_1c_{21}\\
+4c_{13}c_{11}q_1c_{12}c_{23}^2c_{21}^2q_3+c_{13}c_{11}c_{21}^3q_2^2c_{23}c_{12}
-c_{13}q_1^2c_{12}^2c_{23}^2c_{21}^2 +c_{11}^3c_{22}^3q_3^2c_{21}
+c_{13}c_{11}c_{21}^2q_2c_{23}q_1c_{22}c_{12}\\
+4c_{13}^2c_{11}c_{21}c_{22}q_1^2c_{23}^2
-3c_{13}c_{11}^2c_{21}^2q_3^2c_{23}^2+c_{11}^2q_3c_{22}^2c_{12}q_2c_{21}^2-2c_{11}c_{21}^2q_1c_{12}^2c_{23}q_3c_{22},
\end{array}
$$}
%{\small
%$$\begin{array}{l}
%2c_{13}c_{11}^2c_{22}^3q_{(1,0)}c_{21}q_{(0,0)}-3c_{13}^2c_{11}c_{21}^2q_{(1,0)}c_{23}^2q_{(0,0)}
%-2c_{13}c_{11}c_{21}c_{22}^2q_{(0,0)}^2c_{23}c_{12}\\
%+2c_{11}^2c_{21}^2q_{(0,1)}^2c_{23}c_{22}c_{12}
%-c_{13}c_{11}^2c_{22}^4q_{(0,0)}^2
%-c_{13}c_{11}^2c_{22}^2q_{(1,0)}^2c_{21}^2-c_{13}c_{11}^2q_{(0,1)}c_{23}q_{(1,0)}c_{22}c_{21}^2\\
% -c_{13}c_{11}^2q_{(0,1)}c_{23}c_{22}^2q_{(0,0)}c_{21}
%-c_{11}^2q_{(0,1)}c_{22}^3c_{12}q_{(0,0)}c_{21}
%+4c_{13}c_{11}q_{(0,0)}c_{12}c_{23}^2c_{21}^2q_{(0,1)}\\+c_{13}c_{11}c_{21}^3q_{(1,0)}^2c_{23}c_{12}
%-c_{13}q_{(0,0)}^2c_{12}^2c_{23}^2c_{21}^2
%+c_{11}^3c_{22}^3q_{(0,1)}^2c_{21}
%+c_{13}c_{11}c_{21}^2q_{(1,0)}c_{23}q_{(0,0)}c_{22}c_{12}\\
%+4c_{13}^2c_{11}c_{21}c_{22}q_{(0,0)}^2c_{23}^2
%-3c_{13}c_{11}^2c_{21}^2q_{(0,1)}^2c_{23}^2+c_{11}^2q_{(0,1)}c_{22}^2c_{12}q_{(1,0)}c_{21}^2-2c_{11}c_{21}^2q_{(0,0)}c_{12}^2c_{23}q_{(0,1)}c_{22},
%\end{array}
%$$}
which is an irreducible polynomial and does not divide (\ref{cc}).}
\end{example}

\section{Traces and Residues in the torus}\label{4}
In this section, we will compare our denominator formulas for
traces in the torus with those that can be obtained by applying the results given in \cite{CDS}. As in
Section \ref{3}, we will be dealing with a system of $k$ generic
Laurent polynomials $ f_i:=\sum_{a\in \AA_i} c_{ia} \,
\mathbf{t}^a$, where $i=1,\ldots,k,$ $\AA_i$ is a finite subset of
$\ZZ^k$, and the $c_{ia}$'s are indeterminates over $\QQ.$

\par
For a given Laurent polynomial $p$ in
$S:=K[t_1,{t_1}^{-1},\ldots,t_k,{t_k}^{-1}],$ where $K$ was
defined in Section \ref{3}, the global residue of the
differential form
$$\phi_p:=\frac{p}{f_1\ldots f_k}\frac{dt_1}{t_1}\wedge\ldots\wedge\frac{dt_k}{t_k}$$
equals, with our notation,
\begin{equation}\label{adf}
{\rm Residue}^T_f(p):={\rm Trace}\big({\rm
Times}_{\frac{p}{J^T_f}}\big),\end{equation} where $J^T_f$ denotes
the affine toric Jacobian $J^T_f:=\det\left(t_j\frac{\partial
f_i}{\partial t_j}\right)_{1\leq i,j\leq k}.$
\par In \cite{CDS}, formulas for the denominator of the rational expression (\ref{adf})
were proposed. In particular, for any Laurent monomial
$\mathbf{t}^a \in S$, by replacing $p$ with $\mathbf{t}^a J^T_f$
in (\ref{adf}), we get that
$${\rm Trace}({\rm Times}_{\mathbf{t}^a})={\rm Residue}^T_f(\mathbf{t}^a J^T_f),$$
and so, the formulas in \cite{CDS} can be used for computing the
denominator of the trace of ${\rm Times}_{\mathbf{t}^a}.$ Now, we
will compare the denominator formula obtained in this way with
our results in Section \ref{3}.

Assume that, for $i=1,\dots, k$, $\mathcal{A}_i =
\Delta_i \cap \ZZ^k$ for an integral polytope $\Delta_i$ in
$\RR^k$ and, without loss of generality, that
$\mathcal{A}_i\subset (\ZZ_{>0})^k$.

Under these assumptions, the support of each of the polynomials
$t_j\frac{\partial f_i}{\partial t_j}$ $(j=1,\dots, k)$ is
$\mathcal{A}_i$, and therefore, the support of $J_f^T$ is
contained in $\Delta \cap \ZZ^k$, where $\Delta:=\Delta_1 +
\cdots + \Delta_k$. Thus, in order to obtain a denominator for
${\rm Trace}({\rm Times}_{\mathbf{t}^a})$ it suffices to get
denominator formulas for the residues in the torus of those
monomials $\mathbf{t}^m$ with $m\in (\Delta \cap \ZZ^k) + a $.

First, we introduce some notation. For each facet of $\Delta$
with primitive inward normal vector $\omega$, let $a_\omega\in
\ZZ$ be defined as
$$a_\omega:= - \min \{\langle b, \omega \rangle: b\in \Delta\},$$
and for every $m\in \ZZ^k$, let
\begin{eqnarray*}
\mu_\omega^{-}(m)&:=& - \min \{ 0, \langle m, \omega \rangle
+a_\omega -1\},\\
\delta'_\omega(m)&:=& \mu_\omega^{-}(m) [\omega^\perp :
L(\AA_1^{\omega}, \dots, \AA_k^{\omega})].
\end{eqnarray*}

Then, \cite[Theorem 3.2]{CDS} states that $\prod_{\omega} {\rm
Res}_{\AA_1^\omega, \dots, \AA_k^\omega}(f_1^{\omega}, \dots,
f_k^{\omega})^{\delta'_\omega(m)},$ where the product runs over
all the primitive inward normal vectors of facets of $\Delta$, is
a denominator for ${\rm Residue}_f^T(\mathbf{t}^m)$. Therefore, if
$$\delta'_\omega:= \max\{\delta'_\omega(m) : m \in (\Delta \cap
\ZZ^k) + a\},$$ the following polynomial is a denominator for
${\rm Trace}({\rm Times}_{\mathbf{t}^a})={\rm Residue}^T_f(\mathbf{t}^a
J^T_f)$: $$\prod_{\omega} {\rm Res}_{\AA_1^\omega, \dots,
\AA_k^\omega}(f_1^{\omega}, \dots,
f_k^{\omega})^{\delta'_\omega}.$$

Let us estimate the exponents $\delta'_\omega$. For $m\in (\Delta
\cap \ZZ^k) + a $, write $m = m_\Delta + a$ with $m_\Delta \in
\Delta \cap \ZZ^k$. We have that $\langle m, \omega\rangle
+a_\omega - 1 =\langle m_\Delta, \omega\rangle +a_\omega +
\langle a, \omega \rangle - 1\ge \langle a, \omega \rangle -1$.

We will consider two cases separetely:

\begin{itemize}
\item $\langle a, \omega \rangle >0$: for every $ m\in (\Delta
\cap \ZZ^k) + a $, the previous inequality implies that $\langle
m, \omega\rangle +a_\omega - 1 \ge 0$, and so, $\mu_\omega^{-}(m)
= 0 $. Therefore, $\delta'_\omega = 0$.

\item $\langle a, \omega\rangle \le 0$:
taking $m_\Delta \in \Delta \cap \ZZ^k$ so that $\langle m_\Delta,
\omega\rangle +a_\omega = 0$, for $m:= m_\Delta +a$, we get
$\langle m, \omega\rangle +a_\omega - 1 = \langle a, \omega
\rangle -1 \le -1$. Then, $\mu_\omega^{-}(m) =1-\langle a, \omega
\rangle \ge 1$ and therefore, $\delta'_\omega = (1-\langle a,
\omega \rangle)L(\AA_1^{\omega}, \dots, \AA_k^{\omega})$.

\end{itemize}

We deduce the following formula for a denominator of ${\rm
Trace}({\rm Times}_{\mathbf{t}^a})$:
\begin{equation}\label{ellos}
\prod_{\omega: \langle a, \omega \rangle \le 0} {\rm Res}_{\AA_1^\omega, \dots,
\AA_k^\omega}(f_1^{\omega}, \dots, f_k^{\omega})^{(1-\langle a,
\omega \rangle)L(\AA_1^{\omega}, \dots, \AA_k^{\omega})}.
\end{equation}

\bigskip

Finally, we will restate the result in Proposition \ref{prp} in
this context: here $\AA = \{ a \}$, $\AA'= \{ 0 \}$ and $\AA_0=\{
0, a\}$. The first resultant appearing in the factorization
equals $1$. On the other hand, it follows from the definition
that $\mu_\omega =  - \min\{ 0, \langle a, \omega\rangle\}$, and
so, $\delta_\omega = -\langle a, \omega\rangle L(\AA_1^{\omega},
\dots, \AA_k^{\omega})$ if $\langle a, \omega\rangle\le 0$ and
$\delta_\omega =0$ otherwise. We conclude that
\begin{equation}\label{nos}
\prod_{\omega: \langle a, \omega \rangle \le 0} {\rm Res}_{\AA_1^\omega, \dots,
\AA_k^\omega}(f_1^{\omega}, \dots, f_k^{\omega})^{-\langle a,
\omega \rangle L(\AA_1^{\omega}, \dots, \AA_k^{\omega})}
\end{equation} is a
denominator for ${\rm Trace}({\rm Times}_{\mathbf{t}^a})$.

\par By comparing (\ref{nos}) with (\ref{ellos}) we get a slightly improvement on the exponent of each of the factors in the denominator.

\section{Multidimensional Residues in $\CC^n$}\label{5}
In the case where the underlying variety is the zero locus of a
regular sequence of $n$ polynomials in $\CC[x_1,\ldots,x_n]$
without zeroes in the infinity, some results of Jouanolou on
generalized discriminants given in \cite{jou} (see also \cite{GKZ}) will allow us to use our
formulas in order to recover known results about residues
(\cite{elk,CDS1,CDS}).

\subsection{The denominator of the residue}

Let $f:= (f_1,\dots, f_n)$ with $f_i:=\sum_{|a| \le d_i} c_{ia}
x^{a}$, $i=1,\dots, n$, be a generic system of polynomials in
$\KK[x_1,\dots, x_n]$ of respective degrees $d_1,\dots, d_n$. Set
$J_f:= \det \left(\partial{f_i} /
\partial{x_j}\right)$ for the Jacobian determinant of the system.
Let us observe that $\deg(J_f) = \rho:= \sum_{i=1}^n d_i - n$.

For every $\beta \in (\ZZ_{\ge 0})^n$, the global residue
associated with the data $(x^\beta,f)$ can be obtained as
$${\rm Residue}_f(x^\beta) = {\rm Trace}({\rm Times}_{x^\beta/J_f}), $$
where the linear maps in the right hand side of the equation are
defined in the quotient ring $\KK[x_1,\dots, x_n]/\langle f_1,\dots,
f_n\rangle$. Applying (\ref{mtsp}), we obtain the following
expression for the residue as a rational function in the
coefficients of $f$:
$${\rm Residue}_f(x^\beta) =
\frac{\frac{\partial {\rm Res}_{D,d_1,\dots,
d_n}(f_0,f_1,\ldots,f_n)} {\partial c_{0\beta}}|_{f_0\to J_f}}{{\rm
Res}_{D,d_1,\dots, d_n}(J_f,f_1,\ldots,f_n)},$$ where $D:= \max\{
|\beta|, \rho\},\, f_0$ is a generic polynomial of degree $D$, ${\rm
Res}_{D, d_1,\dots, d_n}$ is the classic resultant of a generic
system of $n+1$ polynomials with respective degrees $D, d_1,\dots,
d_n$, and $\frac{\partial}{\partial c_{0\beta}}$ stands for the
derivative with respect to the coefficient of the monomial $x^\beta$
in the first polynomial.

Due to Proposition \ref{prp}, if $f_i^0:= \sum_{|a| = d_i} c_{ia}
x^{a}$ is the homogeneous part of degree $d_i$ of the polynomial
$f_i$, $i=1,\dots, n$,  the denominator of this expression
factors as
\begin{equation}\label{id1}
{\rm Res}_{D, d_1,\dots, d_n}(J_f, f_1,\dots, f_n) = {\rm Res}_{d_1,\dots,
d_n}( f_1^0,\dots, f_n^0)^{D-\rho} \, {\rm Res}_{\rho, d_1,\dots,
d_n}(J_f, f_1,\dots, f_n).
\end{equation}
Now, according to \cite{jou}, the following identity holds:
\begin{equation}\label{jjou}
{\rm Res}_{\rho, d_1,\dots, d_n}(J_f, f_1,\dots, f_n) =
{\rm Res}_{d_1,\dots, d_n}(f_1^0,\dots, f_n^0)\, {\rm Disc}(f).
\end{equation}
Here, ${\rm Disc}(f)$ denotes the discriminant of the polynomial
system $f = (f_1, \dots, f_n)$. Then, we can factor the
denominator of the residue further. This result

\smallskip
\begin{remark}
Identity (\ref{jjou}) can be also proved by applying the results about \textit{principal $A$ determinants}
given in \cite[Chapter $10$]{GKZ} to the polynomial $f_1+\sum_{i=2}^n T_i\,f_i$, with $T_2,\ldots,T_n$ new variables.
Theorem $1.2$ in \cite[Chapter $10$]{GKZ} essentially implies (\ref{jjou}).
\end{remark}

By replacing (\ref{jjou}) in (\ref{id1}), we get
\begin{equation}\label{id2}
{\rm Residue}_f(x^\beta) = \frac{ A_\beta(f)}{{\rm Res}_{d_1,\dots,
d_n}( f_1^0,\dots, f_n^0)^{D+1-\rho}\,{\rm Disc}(f)},
\end{equation}
with
$$A_\beta(f)=\frac{\partial {\rm Res}_{D,d_1,\dots, d_n}(f_0,f_1,\ldots,f_n)}{\partial_{0\beta}}|_{f_0\to J_f}$$
Let us show now that $A_\beta(f)$ is also divisible by the discriminant ${\rm Disc}(f)$.
Recall that ${\rm Disc}(f)$ is an irreducible polynomial in the coefficients of $f$ which vanishes if and only if the
system $f$ has a multiple root  provided that $f$ does not have roots at infinity (see \cite{GKZ}).
\par
Suppose that we specialize $f$ in such a way that ${\rm
Res}_{d_1,\ldots,d_n}(f_1^0,\ldots,f_n^0)\neq0$ but ${\rm Disc}(f)$
vanishes. Hence, the expression on the left hand side of (\ref{id2})
still has sense and is a finite complex number. On the right hand
side of this identity, the denominator vanishes. As ${\rm
Residue}_f(x^\beta)$ is a continuous function of $f$ outside of
${\rm Res}_{d_1,\ldots,d_n}(f_1^0,\ldots,f_n^0)=0$ (see for instance
the ``Principle of Continuity'' in \cite[p.~657]{GH}), then
$A_\beta(f)$ must be zero also. This implies that ${\rm Disc}(f)$
divides $A_\beta(f)$, and hence

\begin{equation}\label{residue}
{\rm Residue}_f(x^\beta) = \frac{ A'_\beta(f)} {{\rm
Res}_{d_1,\dots, d_n}(f_1^0,\dots,f_n^0)^{D-\rho+1} },
\end{equation} with $A'_\beta(f)\in \KK[c_{ia}]$ (c.f.
\cite[Proposition 2.1]{elk}, \cite[Proposition 5.15]{BGVY}).

\begin{remark}
In \cite{CDS}, a formula like (\ref{residue}) is given for systems of generic sparse polynomials,
where the denominator is now a product of sparse \textit{facet} resultants. Our technique is limited only to the
homogenous (dense) case, due to the fact that there is still not known formula like (\ref{jjou}) for sparse systems.
In \cite[Chapter 10]{GKZ} some general formulae is given, but only in the case when the toric variety associated to
the sparse system is smooth.
\end{remark}
\subsection{An algebraic proof of the Euler-Jacobi vanishing theorem}
The following result is very well known in the literature, see
for instance \cite{mac} or \cite[Chapter $5$, Corollary $4$]{AGV}.
We will recover it by using Identity (\ref{residue}) and
algebraic methods.
\begin{theorem}[Euler-Jacobi Formula]
For a system $f=(f_1,\ldots,f_n)$ as before, and any polynomial
$h$ of degree less than the degree of the Jacobian of the system,
$${\rm Residue}_f(h)=0.$$
\end{theorem}

\begin{proof}
It suffices to show that ${\rm Residue}_f(x^\beta) = 0$ for every
monomial of degree at most $\rho:= \deg(J_f)$.

For $i=1,\dots, n$, let $f_{t,\,i}:=\sum_{|a| \le d_i} c_{ia} \,
t^{d_i-|a|} x^{a}\in \KK[t][x_1,\dots, x_n]$ and let $f_t$ be the
generic system $f_t:=(f_{t,1}, \dots, f_{t, n})$.

Fix $\beta \in (\ZZ_{\ge 0})^n$ with $|\beta| <\rho$. In what
follows, we will relate the global residue of $x^\beta$ with
respect to the original system $f$ with its global residue with
respect to the system $f_t$.

First, let us observe that
$$\{ \xi  = (\xi_1,\dots, \xi_n) \in \mathbb{A}^n: f_t(\xi) = 0\} =
\{ t \eta= (t \eta_1,\dots, t \eta_n) \in \mathbb{A}^n: f(\eta) =
0\} $$ and that $J_{f_t}(t\eta) = t^\rho J_f(\eta)$, which implies
that
$${\rm Residue}_{f_t}(x^\beta) = \sum_{f_t(\xi) =
0}\frac{\xi^\beta}{J_{f_t}(\xi)} = \sum_{f(\eta) =
0}\frac{(t\eta)^\beta}{t^\rho J_{f}(\eta)}  = t^{|\beta| - \rho}
\sum_{f(\eta) = 0}\frac{\eta^\beta}{J_{f}(\eta)} = t^{|\beta| -
\rho} \, {\rm Residue}_{f}(x^\beta).
$$
On the other hand, taking into account that $f_i^0 = f_{t,\,
i}^0$, Identity (\ref{residue}) applied to system $f_t$ states
that
$${\rm Residue}_{f_t}(x^\beta) = \frac{ A'_\beta(f_t)} {{\rm
Res}_{d_1,\dots, d_n}(f_1^0,\dots,f_n^0)^{D-\rho+1} }.$$ Finally,
combining the previous identities we deduce:
$$ A'_\beta(f_t) = t^{|\beta| - \rho}  A'_\beta(f).$$

Now, $A'_\beta(f_t)$ is a polynomial in $\KK[t, c_{ia}]$, while
$A'_\beta(f)\in \KK[c_{ia}]$ does not depend on $t$ and $|\beta| -
\rho<0$. Then, the above equality is only possible if $A'_\beta(f)
=0$ and so, ${\rm Residue}_{f}(x^\beta) = 0$.
\end{proof}

\medskip

\noindent {\sc Carlos D'Andrea:}  Universitat de Barcelona,
Departament d'{\`A}lgebra i Geometria. Gran Via 585, 08007
Barcelona, Spain. E-mail: {\tt cdandrea@ub.edu}.

\smallskip
\noindent {\sc Gabriela Jeronimo:} Departamento de
Matem\'atica,
FCEyN, Universidad de Buenos Aires, Ciudad Universitaria, 1428
Buenos Aires, Argentina. E-mail: {\tt jeronimo@dm.uba.ar}.

\end{document}